\documentclass[
11pt,%
tightenlines,%
twoside,%
onecolumn,%
nofloats,%
nobibnotes,%
nofootinbib,%
superscriptaddress,%
noshowpacs,%
centertags]%
{revtex4}
\usepackage{abc}
\begin{document}

\noindent Khrapov S.S., Khoperskov A.V. Application of Graphics Processing Units for Self-Consistent Modelling of Shallow Water Dynamics and Sediment Transport // Lobachevskii Journal of Mathematics, 2020, Vol. 41, No. 8, pp. 1475-1484.

\titlerunning{----} 
\authorrunning{-----} 

\title{Application of Graphics Processing Units for self-consistent modelling of shallow water dynamics and sediment transport}

\author{\firstname{S.~S.}~\surname{Khrapov}}
\email[E-mail: ]{khrapov@volsu.ru}
\affiliation{Volgograd State University, Volgograd, Russian Federation}

\author{\firstname{A. V.}~\surname{Khoperskov}}
\email[E-mail: ]{khoperskov@volsu.ru}
\affiliation{Volgograd State University, Volgograd, Russian Federation}


\firstcollaboration{(Submitted by A.~A.~Editor-name)} 


\begin{abstract} 
In this paper, we describe a numerical algorithm for the self-consistent simulations of surface water and sediment dynamics. The method is based on the original Lagrangian-Eulerian CSPH-TVD approach for solving the Saint-Venant and Exner equations, taking into account the physical factors essential for the understanding of the shallow water and surface soil layer motions, including complex terrain structure and its evolution due to sediment transport. Additional Exner equation for sediment transport has been used for the numerical CSPH-TVD scheme stability criteria definition. By using OpenMP-CUDA and GPUDirect technologies for hybrid computing systems (supercomputers) with several graphic coprocessors (GPUs) interacting with each other via the PCI-E / NVLINK interface we also develop a parallel numerical algorithm for the CSPH-TVD method. The developed parallel version of the algorithm demonstrates high efficiency for various configurations of Nvidia Tesla CPU + GPU computing systems. In particular, maximal speed up is 1800 for a system with four C2070 GPUs compare to the serial version for the CPU. The calculation time on the GPU V100~(Volta architecture) is reduced by 95 times compared to the GPU C2070~(Fermi architecture).  
\end{abstract}



\maketitle


\section{Introduction}

Modelling of unsteady flows of the surface fluid requires the calculation of movable flow boundaries \cite{Bulatov-Elizarova2016} where the presence of sediment transport complicates this problem because the boundaries between the fluid and sediment may not coincide. Another complicating factor, in addition to the heterogeneity of the bottom topography $b(x, y)$, is the heterogeneity of soil characteristics (e.g., porosity, particle sizes). A detailed analysis of various models of sediment transport can be found in \cite{Parker2008}.  Usually, the mathematical model of a shallow water dynamics is based on hyperbolic equations for which Godunov-type solvers are often used \cite{Kulikovskii-etalBook}.  The first family of solvers uses the Harten-Lax-Leer~(HLL) scheme and its modifications \cite{Harten-Hyman1983, Toro2019}. Another approach is based on Roe's approximate solver~\cite{Ambrosi1995}. Nowadays, another state-of-the-arm approach for modelling sediment transport by water flow is the discrete least-squares meshless method \cite{Erami-Firoozjaee2020}.

Recently advantages of the architecture of modern GPUs have been successfully used to solve a wide range of mathematical problems~\cite{Becerra-Sagredo-etal2020, Chisholm-etal2020, Erofeev-etal2019, Shaikhislamov-etal2019, Dazzi-etal2019, Lastra-etal2018, Bazulin-etal2018} thus making usage of the parallel computing on GPUs even more topical for a broad range of academical and applied problems. In this work, we discuss a combined SPH-TVD method for the numerical integration of the Saint-Venant and Exner equations which was proposed previously in our works \cite{Khrapov-etal2013, Khoperskov-Khrapov2018}. In this paper, we extend our approach for the modelling of self-consistent dynamics of surface water and sediment.  In particular, the main aim of our work is to develop an effective parallel code for various computing systems with GPUs.

\section{Математическая модель}
It is widely accepted that a consistent description of the dynamics of the surface water flow and sediment transport requires the use of multilayer models in the shallow water approximation \cite{Broecker-etal2019, Sukhinov-etal2018, Elizarova-Ivanov2018}.
In this approach, the Saint-Venant equations should be written in the form of the Exner equation \cite{AlanPeter2008, Dazzi-etal2019, Erami-Firoozjaee2020} which allows simulating the dynamics of soil erosion and sediment movement due to fluid flow. Deformations of the bottom surface in weak currents are accompanied by the formation of structures in the form of ripples, ridges, dunes, etc., which are an important structural element of the bottom, having a significant effect on hydraulic flow resistance \cite{DyakonovaKhoperskov2018}. For example, strong currents lead to a significant erosion of the bottom, in particular, a channel can be formed due to the dam break and then rapidly deepens and expands under the influence of increasing of the fluid flow.

The Exner equation for one-component soil model has the form:
\begin{equation}\label{eq:erosion}
(1-\psi)\,\frac{\partial b}{\partial t} + \nabla_{\perp} {\bf J}_b = q_{b}^{+} - q_{b}^{-}\,,
\end{equation}
where ${\bf J}_b$ is the vector of solid transport discharge, $q_{b}^{+}$ is the matter source due to sediment deposition, $q_b^{-}$ is the drain function due to bottom erosion, $\psi$ is the bottom porosity, $\nabla_{\perp} = \left\{ \partial/\partial x, \partial/\partial y \right\}$ is the component of the operator nabla in the Earth plane, $(x,y)\perp z$.

 Complete modelling of the erosion processes, soil transport and sedimentation is a challenging mathematical and computational problem \cite{Ipatov-etal2019, Anuprienko-Kapyrin2018} and its solution is only possible by using the phenomenological approach complemented by a careful preliminary examination of the simulation results and the benchmark experiments \cite{Coronel-etal2019, Torok-etal2019}. 
The equation (\ref{eq:erosion}) is nonlinear due to the variety of physical processes in the bottom layer.
   For example, the flow of material ${\bf J}_b$ can depend both on the force of the liquid friction against the bottom and the slope of the bottom \cite{Liu2008}:
\begin{equation}\label{eq:erosionGb}
	\mathbf{J}_b =  \mathbf{J}^0_b - C_J\,{\rm sign}(\mathbf{J}^0_b) |\mathbf{J}^0_b| {\nabla}_\perp b\,,
\end{equation}
where $\mathbf{J}^0_b$ is the solid transport discharge for zero inclination, $C_J = 1.5\div 2.3$ is an arbitrary empirical constant that can reach values  $C_J=5$ \cite{Jayaratne-etal2012}. 
Thus, the equation (\ref{eq:erosion}) contains a nonlinear diffusion term in the form of $\displaystyle {\nabla}_\perp |\mathbf{J}^0_b| {\nabla}_\perp b$ taking into account the relation (\ref{eq:erosionGb}).
To determine the flux $ \mathbf{J}^0_b$, we use the Grass formula \cite{Liu2008}:
\begin{equation}\label{eq:Grass}
\mathbf{J}^0_b = A_J\, \mathbf{v}\,|\vec{v}|^m \,,
\end{equation}
where $\mathbf{v}$ is the average lateral flow velocity above the bottom, coinciding with the fluid velocity in the shallow water model, $A_J, m$ are the constant coefficients, for example, $A_G = 0.001, m = 2$ for fine sand. 

  The value of $A_J$ in the ratio (\ref{eq:Grass}) which depends on the characteristics of the moving substance and it can be introduced as follows \cite{Debolskaya-etal2013}:
\begin{equation}\label{eq:AG}
A_J = \frac{0.05\, n_M^3}{(s-1)\sqrt{g H} \,d_{50}}\,,
\end{equation}
where $s=\rho_s/\rho$ is the relative mass density, $d_{50}$ is the median sediment particle sizes.
 The erosion of the bottom soil occurs when the velocity exceeds the critical value determined by the Shamov's equation~\cite{Przedwojski-etal1995}:
\begin{equation}\label{eq:Ucr}
v_k = C_{Sh}\, d_{50}^{1/3}\, H^{1/6}\,, \quad C_{Sh} \simeq 4\div 6.
\end{equation}

The vector form of the system of shallow water equations and Exner equation can be written as follows:
\begin{equation}\label{eq:vectorform}
\frac{\partial {\bf U}}{\partial t} +  \frac{\partial {\bf F}}{\partial x} + \frac{\partial {\bf G}}{\partial y} = {\bf \Phi} \,,
\end{equation}
where
$$
\mathbf{U} = 
\begin{pmatrix}
H    \\ 
Hu   \\ 
Hv   \\
b
\end{pmatrix}, \quad
\mathbf{F} = 
\begin{pmatrix}
Hu    \\ 
Hu^2   \\ 
Huv   \\
\displaystyle\frac{J_{b x}}{1-\psi}
\end{pmatrix},\quad
\mathbf{G} = 
\begin{pmatrix}
Hv    \\ 
Huv   \\ 
Hv^2   \\
\displaystyle\frac{J_{b y}}{1-\psi}
\end{pmatrix}, \quad
{\bf \Phi} = 
\begin{pmatrix}
\sigma    \\ 
\displaystyle gH \nabla_\perp (H+b) + H \mathbf{f}   \\ 
\displaystyle\frac{q_b^+ - q_b^-}{1-\psi}
\end{pmatrix},
$$
where $H$ is water depth, $u(x,y,t),\,v(x,y,t)$ is the water velocity components, $\sigma$ is the function of water sources/drains, determined by the operating mode of hydraulic structures, rainfall/snowfall, the rainwater runoff, the infiltration/evaporation; ${\bf f}$ is the external forces connected with the Earth rotation, the wind, the bottom friction, the turbulence.

We use the Combined Smoothed Particle Hydrodynamic method~(CSPH-TVD), which contains the Lagrangian stage~(modified Smoothed Particle Hydrodynamic, or SPH) and the Euler stage (modified Total Variation Diminishing, or TVD)~\cite{Khrapov-etal2013}. The consistent use of these two stages during the integration provides several important positive properties that can be achieved with less computational resources \cite{Dyakonova-etal2016, Khoperskov-Khrapov2018}.
The equation~(\ref{eq:erosion}) imposes additional constraints on the stability of the entire numerical CSPH-TVD algorithm compared to the single-layer shallow water model considered in \cite{Dyakonova-etal2016, Khoperskov-Khrapov2018}. Nonlinear terms of the diffusion type in the equation (\ref{eq:erosion}) due to the relation (\ref{eq:erosionGb}) affect the stability condition, which takes the form:
\begin{equation}\label{Eq2_stability}
\tau_n  =   K \,    \min\limits_{\mathbf{i}} \left(  \frac{h}{2\,\upsilon_p} \, ,     \frac{h}{\upsilon_s}\,,  \frac{h^2}{2D} \right) \, ,
\end{equation}
where $0<K<1$ is the Courant number, $h$ is cell size, $\mathbf{i}$ is the vector spatial index, $\upsilon_p =   \max\limits_\mathbf{i} \left| \mathbf{v}_\mathbf{i}\right|$,  
$\upsilon_s =   \max\limits_\mathbf{i} \left( \left| \mathbf{v}_\mathbf{i}\right| + \sqrt{gH_\mathbf{i}} \right)$, 
$\displaystyle D = \max\limits_\mathbf{i} \left| \mathbf{J}_\mathbf{i}^0\right|/(1-\psi)$.

\section{Parallel implementation of numerical model}

\begin{figure}[!t]
	\centering\includegraphics[width=0.8\hsize]{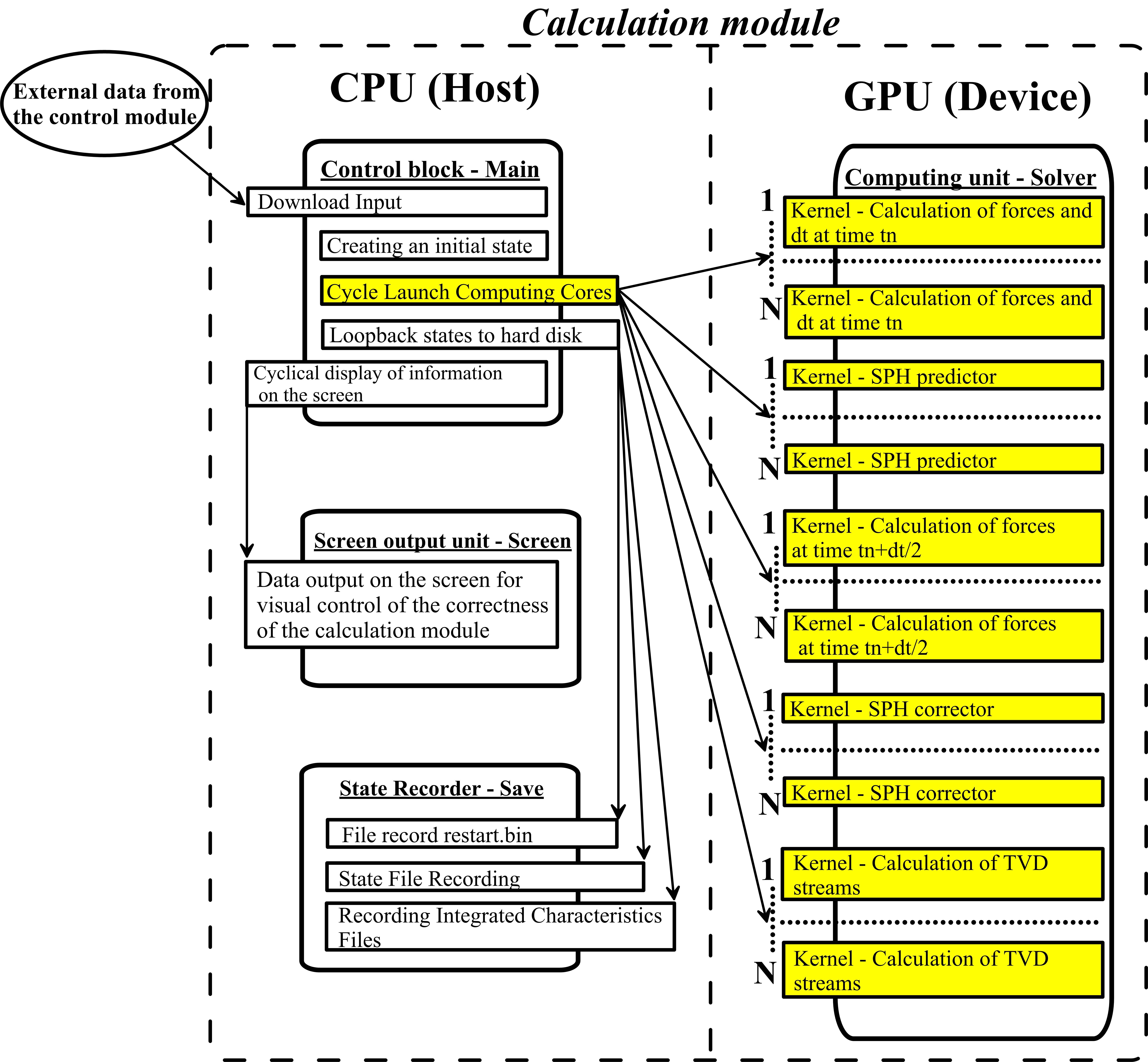}
	\caption{Components diagram of computing module where $N$ is the number of parallel copies of the kernels launched on a GPU device.}
	\label{fig:SW_CUDA}
\end{figure}

 The computational module of the software package for simulation of dynamics of the water surface and sediment based on the CSPH-TVD method consists in four main blocks: Control Block --- $Main$; Computing Block --- $Solver$; State Recording Block --- $Save$; Screen Output Block --- $Screen$ (Fig.~\ref{fig:SW_CUDA}). The Control Block $Main$ provides loading of input data (parameters of the numerical model, the digital elevation model $b(x,y)$, the distributions of bottom friction coefficient $n_M(x,y)$ and absorption coefficient $\beta(x,y)$) and creating of the initial state.  This unit is stored on the CPU, and the computation is organized as follows: the cyclical launch of the computing cores from the $Solver$ block located on the GPU; the cyclical recording of calculation states by calling routines from the block $Save$ (CPU); the cyclical display of the data necessary for monitoring the correct operation of the application by accessing the block $Screen$ (CPU).

The $Solver$ block contains $5$ main computing kernels organized in the form of subprograms and functions where the CSPH-TVD numerical algorithm for computation the shallow water dynamics on irregular terrain and the sediment transport is written in CUDA C. Since the stages~(kernels) of the task should be performed sequentially~(i.e. the data obtained at each individual stage are used in the subsequent one), we chose a parallelization model for the algorithm implementation. Implementation of this approach on GPUs allows us to speed-up the program by hundreds and thousands of times, depending on the configuration of the computing cluster in comparison to the CPU. The task parallelization model (CUDA streams) can be used to separate computation processes and copy data between the CPU and GPU. Since the implementation of the algorithm focuses on the minimization of the traffic between the CPU and the GPU~(data exchange between the CPU and GPU occurs only after 100--1000 time iterations), the parallelization can increase productivity by 3 -- 10\% at best. The speed of the data transfer process between the CPU and the GPU is an order of magnitude lower than the computations time on the GPU. The number of simultaneously running computing kernels for each stage of the SPH-TVD method in the case of the parallel OpenMP version~(for the CPU) is equal to the number of processor kernels $k$ (for systems with shared memory $k = 2,4, ..., 16$). At the same time, the CPU load of our task is about 90~\%, and the speed-up of the computations is $\sim (0.7-0.8)\, k$.

   The number of simultaneously running kernels for each stage of the SPH-TVD method for a CUDA program on a GPU does not depend on the number of scalar GPU kernels (Tesla C1060 --- 240 kernels, Tesla C2070 --- 448 kernels, Tesla K20 --- 2496 kernels, Tesla K40 --- 2880 kernels, Tesla K80 --- 4992 kernels, Tesla P100 --- 3584 kernels, Tesla V100 --- 5120 kernels) and it can be set based on the principle that ``the more -- the better'' the loading of the GPU. In our case, this number is equal to the number of cells of the computational grid, i.e. each thread of the GPU processes a single element. The GPU load of our task is about 99.9\%, and the speed-up of the program~(compared to the serial version of the CPU) reaches $\sim \alpha k$ times ($k$ is the number of scalar GPU cores, $\alpha$ is the ratio of CPU frequency to GPU frequency). The efficient usage of the shared and global memory, the data alignment, and the vectorized data types for the Fermi, Kepler, Pascal, and Volta architectures allow speeding up the program in $\sim (2.5 \div 3) \, \alpha k$ times, respectively.

\begin{figure}[!t]
	\centering\includegraphics[width=0.7\hsize]{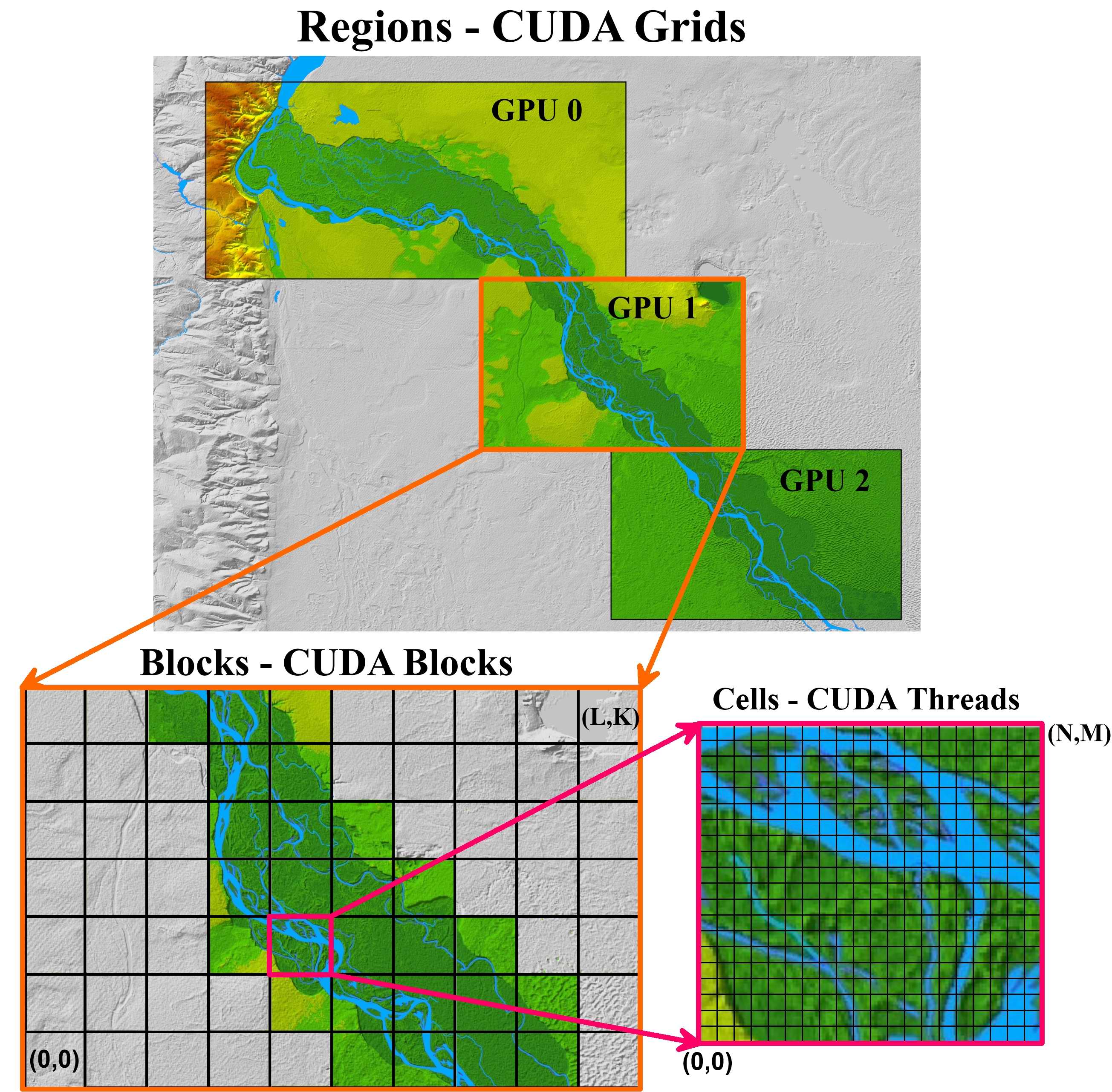}
	\caption{ 
	Diagram demonstrates the structure of the HGS and its connection to the CUDA program model~\cite{Dyakonova-etal2016}.}
	\label{fig:ISS_CUDA}
\end{figure}

Our program has been developed in CUDA C language using the CUDA Runtime libraries. To improve the performance of the numerical algorithm we also adopted a hierarchical grid system (HGS). The HGS-based approach implies that the cells of the computational grid of size $N \times M$ are grouped into blocks of size $L \times K$ which in turn are joint into regions (Fig.~\ref{fig:ISS_CUDA}). Similar to the computation module for the CPU, the application of the HGS approach to the CUDA program leads to a significant speed-up of the computations, especially in case of large sections area without water. Execution of computational kernels in the CUDA version of the code is organized as a grid of thread blocks. Figure~\ref{fig:ISS_CUDA} demonstrates that stream blocks in CUDA are analogous to blocks in HGS. Thereby to implement the HGS algorithm in the GPU version of the program, it is enough to enter the array Index\_block [N\_block] (where N\_block is the number of CUDA blocks in the CUDA grid), which determines the number of cells without water in the CUDA block at a given time, and to supplement the algorithm by new CUDA-kernel, in which the values of the elements of this array are derived.
The body of the corresponding computational core is launched and fully loads the streaming multiprocessor (SM) of the GPU if the condition Index\_block [k\_block]~$>0$ is satisfied (where k\_block = blockIdx.x + blockIdx.y * gridDim.x is the current number of the CUDA block, determined through the built-in CUDA Runtime variables), otherwise SM is freed to start another CUDA kernel from the CUDA Runtime queue.
Code snippets for computing kernels using HGS are shown below.

\medskip
\noindent
{\it CUDA-kernels code for Lagrangian stage}
\begin{verbatim}
__global__ void kernel_LS(double *H, double *Zb, double2 *HV, 
int *Index_block,...){
if(Index_block[blockIdx.x+blockIdx.y*gridDim.x] > 0){
[----//---- COMPUTING KERNEL BODY] } }
\end{verbatim}

\medskip
\noindent
{\it CUDA-kernels code for the Euler stage}
\begin{verbatim}
__global__ void kernel_ES(double *H, double *Zb, double2 *HV, 
int *Index_block,...){
int ib=blockIdx.x+blockIdx.y*gridDim.x;
int xbl=blockIdx.x-1;  if(xbl<0){ xbl=0;}
int ibl=xbl+blockIdx.y*gridDim.x;
int xbr=blockIdx.x+1; 
if(xbr>gridDim.x){xbr=gridDim.x;}
int ibr=xbr+blockIdx.y*gridDim.x;
int ybl=blockIdx.y-1;  if(ybl<0){ybl=0;}
int ibb=blockIdx.x+ybl*gridDim.x;
int ybr=blockIdx.y+1;
if(ybr>gridDim.y){ybr=gridDim.y;}
int ibt=blockIdx.x+ybr*gridDim.x;
if(Index_block[ib] > 0 || Index_block[ibl] > 0 
|| Index_block[ibr] > 0 || Index_block[ibb] > 0 
|| Index_block[ibt] > 0){
[----//---- COMPUTING KERNEL BODY] } }
\end{verbatim}

\medskip
\noindent
{\it CUDA-kernel code Index\_block}
\begin{verbatim}
__global__ void kernel_Index_block(double *H, int *Index_block,...){
if(H[ind]>Eps || Index_Q[ind]>0) Sij[ind_thb]=1;
else Sij[ind_thb]=0;   __syncthreads();
int k = ithbx*ithby/2;
while(k != 0){
if(ind_thb < k) Sij[ind_thb] += Sij[ind_thb+k];
__syncthreads();  k /= 2; }
if(ind_thb == 0) Index_block[ib] = Sij[0];} }
\end{verbatim}

 It should be noted that at the Lagrangian stage it is enough to check the presence of water only in the current CUDA block. Meanwhile, the Euler stage requires the presence of water in neighbouring CUDA blocks to be checked well because fluxes of mass and momentum are calculated at the cell boundaries and water can flow from one block to another. Solution of any applied problems with a large number of calculations is indeed requires the hybrid computing systems with GPUs. For instance, MPI technology provides parallelization of applications on computing systems with shared memory. From other hand, for systems with global memory based on a hybrid configuration of CPU + $k\times$GPU~(along with CUDA technology) it is convenient to use OpenMP because it is necessary to create a separate threads on the CPUs to load each GPU. We used NVIDIA GPUDirect (Peer-to-Peer Communication) technology to exchange data between GPUs without copying them into the CPU when implementing a program adapted for computing on hybrid supercomputers (CPU + $k \times$GPU) using a two-level parallelization scheme OpenMP--CUDA. This approach leads to an almost linear speed-up with an increase in the number of GPUs used. The following is an example of a code fragment that is based on the parallel computing technologies mentioned above.

\medskip
\noindent
{\it Sample Code OpenMP--CUDA}
\begin{verbatim}
#pragma omp parallel for schedule(static,1) shared(deviceCount,
threads,dt)
private(deviceId,devId_b,devId_t)
for(deviceId=0; deviceId<deviceCount; deviceId++){
cudaSetDevice(deviceId);
dim3 blocks(Nx/ithbx,Ny_dev[deviceId]/ithby);
if(deviceId==0){devId_t=1; devId_b=deviceCount-1;
}else if(deviceId==deviceCount-1){devId_t=0; devId_b=deviceId-1;
}else{devId_t=deviceId+1; devId_b=deviceId-1;}
cudaDeviceEnablePeerAccess(devId_t, 0);
cudaDeviceEnablePeerAccess(devId_b, 0);
kernel_forces_predictor<<<blocks, threads>>>( dev_outSrc[deviceId],
dev_inSrc[deviceId],dev_inSrc[devId_b],dev_inSrc[devId_t],
dt,deviceId,deviceCount); }
\end{verbatim}

We increase the total number of CUDA kernels of the calculation module to $8$~(see Fig.~\ref{fig:SW_CUDA}) by adding a new CUDA kernel to implement the HGS approach which
includes the  calculations of the time step and the updating the system at the Eulerian stage in separate CUDA kernels. This extra module allows us to use more efficiently the resources of the graphics processor when modeling the joint dynamics of shallow water and sediment transport. Below we list the CUDA kernels of the calculation module:
  
\noindent K1 determines the presence of water in the CUDA block;

\noindent K2 calculates the forces at the Lagrangian stage of the CSPH-TVD method at time $t_n$ (predictor);

\noindent K3 calculates the time step $dt_{n}$;

\noindent K4 determines the positions of particles and their integral characteristics on the time layer $t_{n + 1/2}$ (predictor);

\noindent K5 calculates the forces at the Lagrangian stage of the CSPH-TVD method at the time $t_{n + 1/2}$ (corrector);

\noindent K6 determines the positions of the particles and their integral characteristics for the next time layer $t_{n + 1}$ (corrector);

\noindent K7 calculates the fluxes of physical quantities through the cell boundaries at the Eulerian stage in $t_{n + 1/2}$;

\noindent K8 updates the state of the system at the Eulerian stage at time instant $t_{n + 1}$, calculating the flow parameters in the cells of the computational grid.
 
\begin{figure}[!t]
	\centering\includegraphics[width=0.55\hsize]{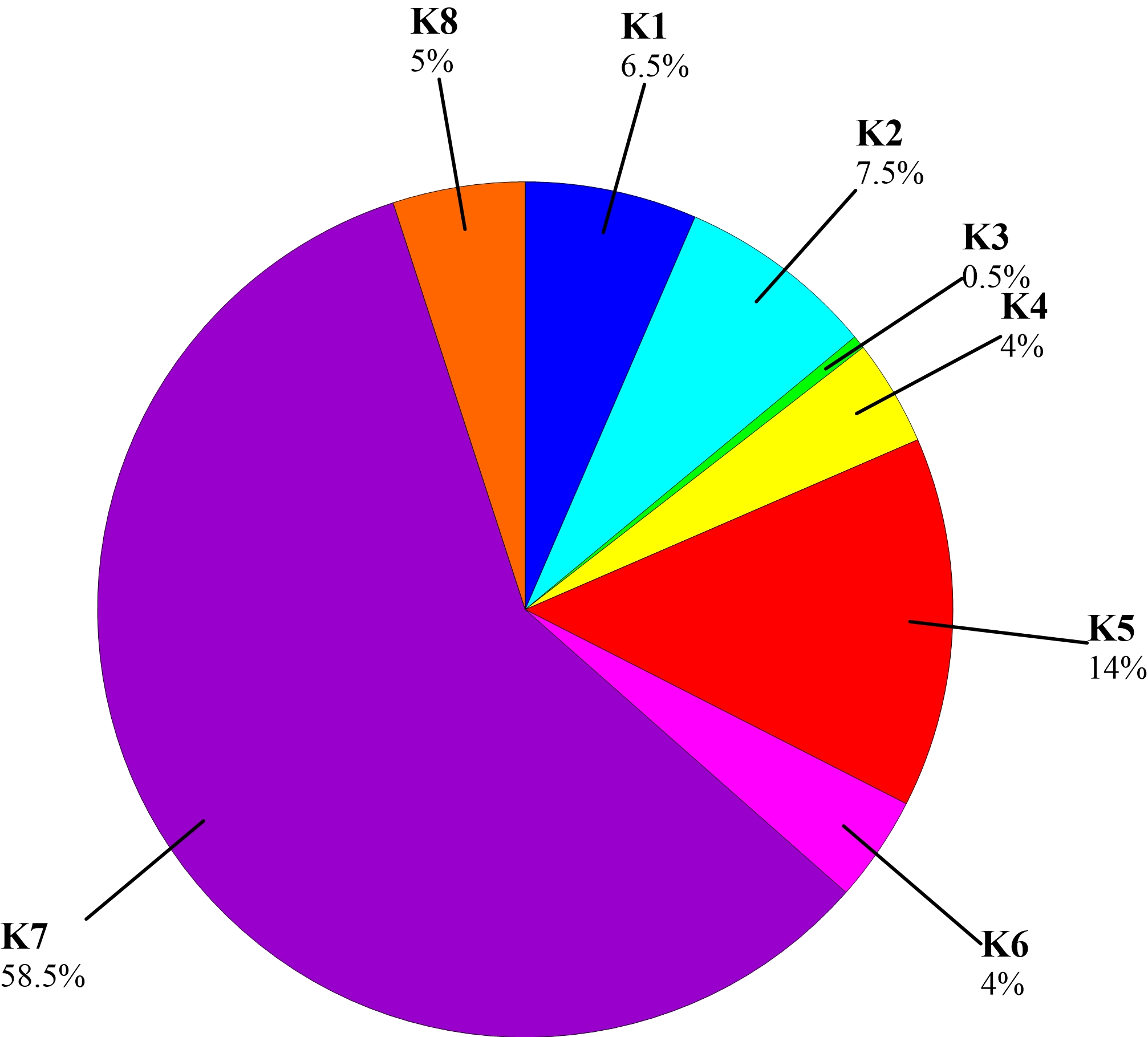}
	\caption{
	Pie chart of the distribution of computational costs on CUDA cores.}
	\label{fig:K1-K8}
\end{figure}

The parallelization efficiency of the CUDA-algorithm of the CSPH-TVD method was studied by simulation the dynamics of seasonal and catastrophic flooding of the Volga-Akhtuba floodplain~(VAF).
The calculations were performed on grids $1024 \times 1024$, $2048 \times 2048$ and $4096 \times 4096$ by using Nvidia Tesla GPUs: C2070, K20, K40, K80 (2~GPU), P100 and V100.

\begin{figure}[!t]
	\centering\includegraphics[width=0.9\hsize]{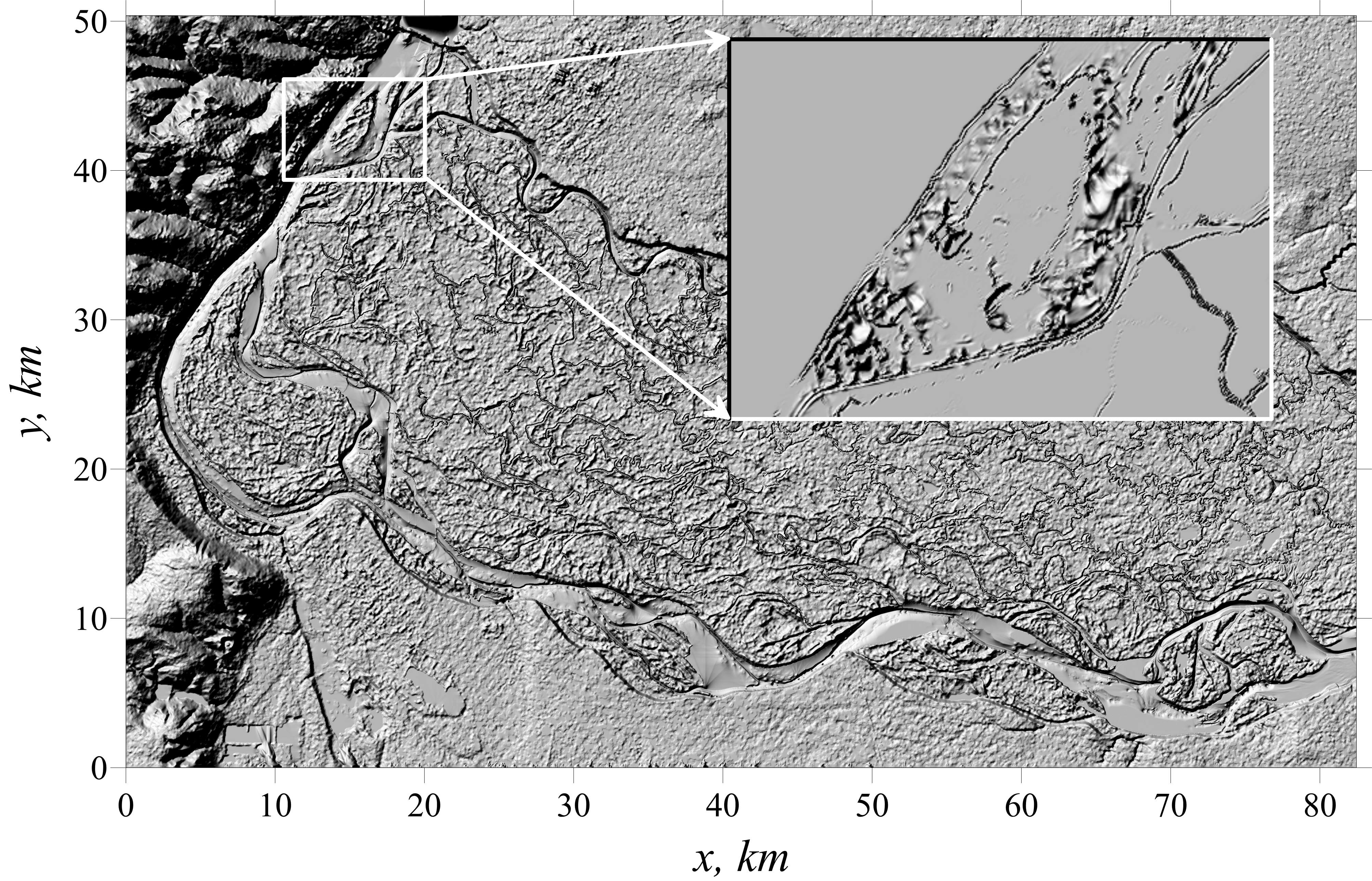}
	\centering\includegraphics[width=0.9\hsize]{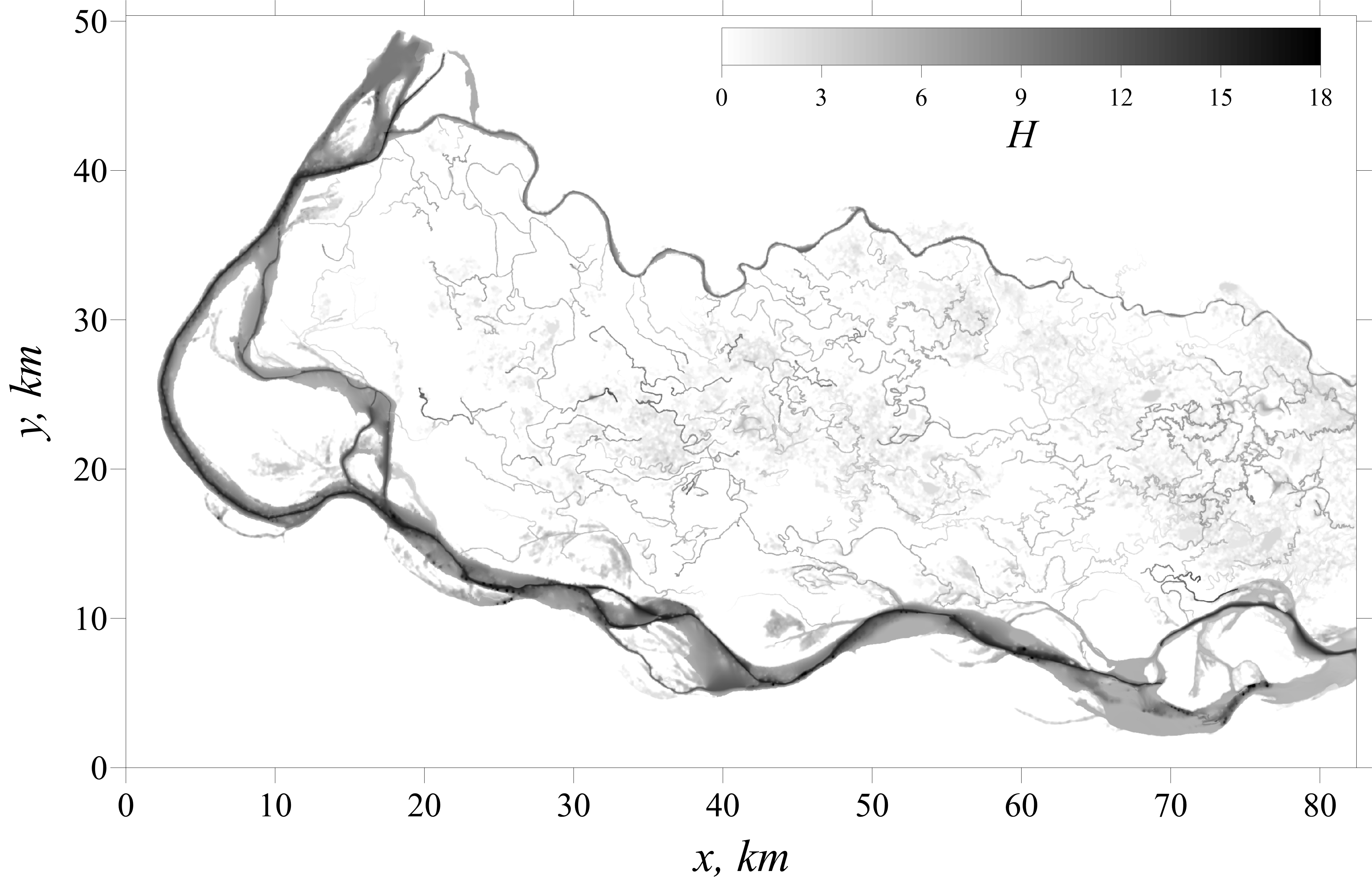}
	\caption{
	Example of modeling the joint dynamics of shallow water and sediment transport in the Volga-Akhtuba floodplain. 
	Top --- Shaded Relief Map, The inset is the change in bottom topography for the 28th day from the beginning of the spring flood. There are the bottom changes within $\Delta b = -1.5 \div 1.5$\,m).
Below is the distribution of water in the floodplain at the same time ([$H$]=m).}
	\label{fig:VAF}
\end{figure}

\begin{figure}[!t]
	\centering\includegraphics[width=0.7\hsize]{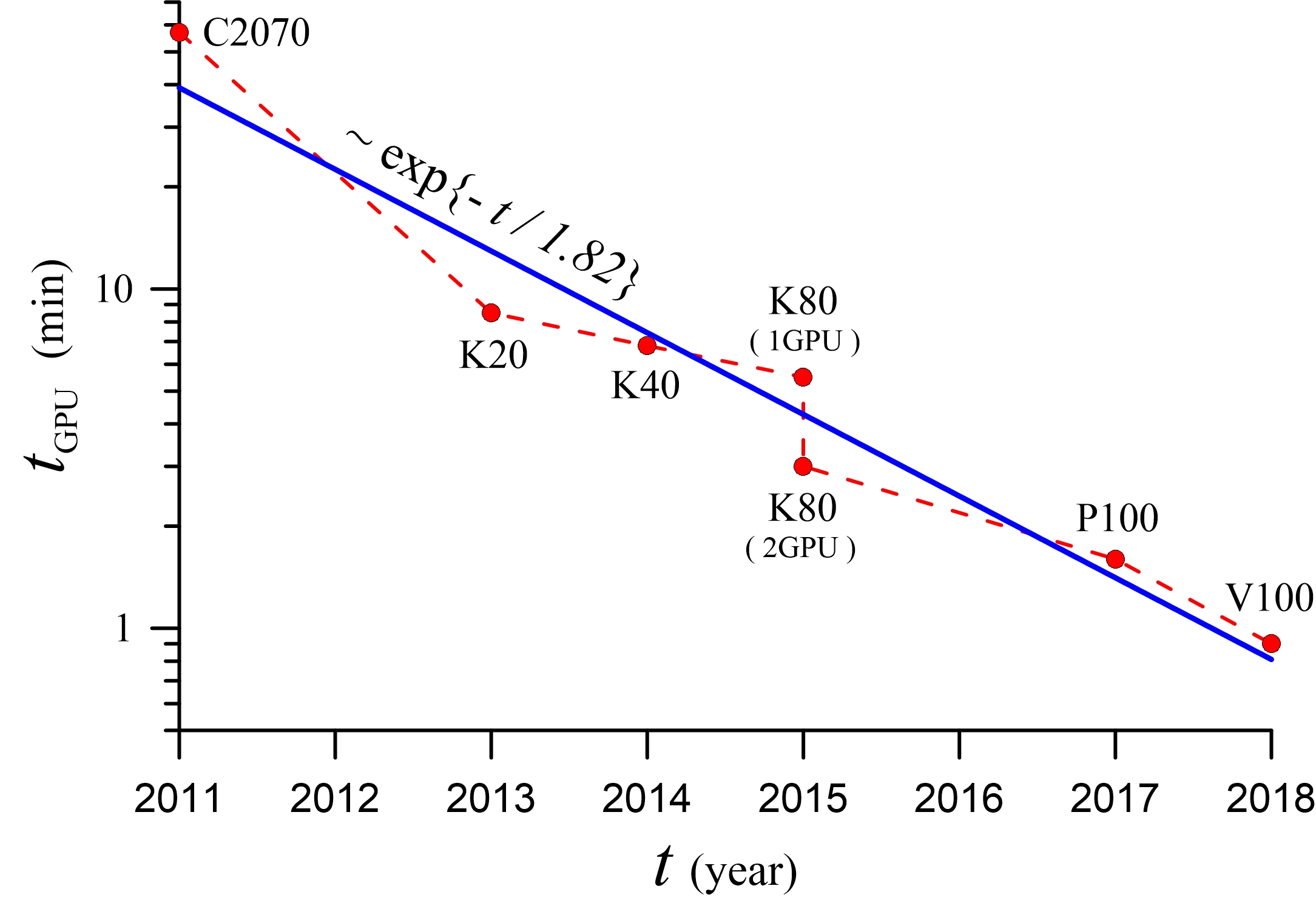}
	\caption{
	Execution time of the parallel CUDA algorithm versus the release year of a new GPU. The solid blue line shows the exponential fit of the relation in the form $\sim \exp \{- t / \tau \}$.
}
	\label{fig:Tgpu(GPU)}
\end{figure}

  Fig.~\ref{fig:K1-K8} shows the distribution of the execution time of the numerical algorithm between CUDA kernels. The most computationally expensive part is CUDA-kernel 7~($58.5\%$) where the fluxes of physical quantities at the boundaries of the Euler cells are calculated. 
This is due to the fact that for each of the x- and y-directions, a piecewise linear reconstruction of the grid functions is carried out and the Riemann problem is solved using the HLL method.
In addition, when calculating the slopes of piecewise-linear distributions of grid functions, TVD-limiters are used, which, on the one hand, suppress oscillations in the numerical solution, and, on the other hand, increase the computation time due to the use of conditional operators in the CUDA software implementation.
  The figure shows that the fraction of the time required for the Lagrangian stage (CUDA kernels K2 -- K6) is $30\%$, the Euler stage~(CUDA kernels K7 -- K8) accounts for $63.5\%$, and the CUDA core K1 that implements the HGS approach, uses only 6.5\%, while it can increase the calculation speed by a factor of ten and more. Thereby in order to increase the computational performance of the algorithm, it is necessary firstly to optimize the kernel K7.  
  Fig.~\ref{fig:VAF} demonstrates the results of simulations of strong spring flooding of the Volga-Akhtuba interfluve, showing changes in the terrain structure due to sediment transport, which are consistent with the analysis of long-term observation series for VAF \cite{Isaeva-etal2019}.
  Fig.~\ref{fig:Tgpu(GPU)} shows a comparison of the computational efficiency of our parallel algorithm of the CSPH-TVD method for various GPUs.  We see that the computation time decreases exponentially ($\sim \exp \{- t / \tau \}$) with the transition to more modern GPUs.
For example, the calculation time on the GPU V100~(Volta architecture) is reduced by 95 times compared to the GPU C2070~(Fermi architecture). We have estimated the characteristic time scale $\tau \simeq 1.82$ for our CUDA-algorithm on the grid $1024 \times 1024$. The parameter $\tau$ decreases with increasing grid resolution (the number of cells) because of the more effiecient loading of new generations of GPUs increases.

We have tested two parallel CUDA versions of the calculation module. The first version can be used on single GPU, in particular on any NVIDIA GPUs with the level of computing capabilities of at least~1.3 (for example, Tesla C1060). The second version is designed for calculations on NVIDIA Tesla C2070, K20, K40, K80, P100 and V100 GPUs with the level of computing capabilities of at least~2.0. Since this version is based on the OpenMP, CUDA and GPUDirect technologies, in principle,  it can utilize all the GPUs and CPU cores of the supercomputer. A particular attention we paid to the effective work with global and shared GPU memory, as well as to the optimization of the data exchange process between the CPUs and GPUs for both versions of our code. For example, the maximum speed-up is 1800 on a supercomputer with four C2070 GPUs compare to the serial version for the CPU and 200 compared to the parallel OpenMP version for the CPU.


\section{Conclusion}

%
%

In conclusion, we summarize the main results of the work.
We have described a new stable numerical algorithm for simulation of self-consistent dynamics of the surface water and sediment process. The method is based on the CSPH-TVD approach by taking into account the process of fluid and grains transport. The implemented parallel CUDA-algorithm CSPH-TVD has also been parallelized by means of direct GPU technology for high-performance supercomputers with GPUs, where the interconnection between nodes is based on PCI-E / NVLINK standards. We demonstrate the high computational efficiency of our algorithm for several different GPU Nvidia Tesla hardware (C2070, K20 / 40/80, P100, V100).

The code described above has been developed for the fast high-resolution simulation of the dynamics of flooding of vast territories. Apart from the calculation module, based on the parallel CSPA-algorithm CSPH-TVD method, our software also includes several routines for the initial conditions preprocessing and visualization of the simulations results. By using this software package, we made a set of new dynamical models of the hydrological regime of the Volga-Akhtuba floodplain were we also take into account sediment transport. We demonstrate the deepening of the bottom of the Volga channel, which, in turn, negatively affects the efficiency of the floodplain territories flooding during the spring floods.

\begin{acknowledgments}
The authors were supported by the Ministry of Science and Higher
Education of the Russian Federation (government task No. 0633-2020-0003). 
 Numerical simulations were carried by using the equipment of the shared research facilities of HPC computing resources at Lomonosov Moscow State University supported by the project RFMEFI62117X0011.
\end{acknowledgments}


%
%
%
%

%
%

\end{document}